# THE TUTTE DICHROMATE AND WHITNEY HOMOLOGY OF MATROIDS

## DAVID G. WAGNER

ABSTRACT. We consider a specialization $Y_M(q,t)$ of the Tutte polynomial of a matroid $M$ which is inspired by analogy with the Potts model from statistical mechanics. The only information lost in this specialization is the number of loops of $M$. We show that the coefficients of $Y_M(1-p,t)$ are very simply related to the ranks of the Whitney homology groups of the opposite partial orders of the independent set complexes of the duals of the truncations of $M$. In particular, we obtain a new homological interpretation for the coefficients of the characteristic polynomial of a matroid.

## 0. INTRODUCTION.

In 1954, Tutte [30] introduced the *dichromate* of a (finite) graph, which has since become known as the *Tutte polynomial*. In the four decades since then this has provided a profound link between combinatorics and other branches of mathematics as diverse as statistical mechanics [3, 5, 17, 19, 20], low-dimensional topology [17, 18, 20], and the theory of Grothendieck rings [9, 10, 29]. Indeed, in 1947 Tutte [29] showed that the "Tutte-Grothendieck ring" $K_0(\mathfrak{G})$ of a suitably defined category $\mathfrak{G}$ of graphs is $\mathbb{Z}[x,y]$; the class $T_G(x,y)$ in $K_0(\mathfrak{G})$ of a graph $G$ is its Tutte polynomial. This construction was axiomatized for "bidecomposition categories" and applied to a category $\mathfrak{M}$ of matroids by Brylawski [9, 10], with the result that $K_0(\mathfrak{M}) = \mathbb{Z}[x,y]$ as well. (It is interesting to compare Brylawski's axiomatization with the usual hypotheses of algebraic $K$-theory; see, *e.g.* Chapter 5 of Silvester [27].) Crapo's generalization [14] of Tutte polynomials to matroids rests on this foundation. The fact that $T_M(x,y)$ is the "universal Tutte-Grothendieck invariant" of the matroid $M$ is just a restatement of the fact that it is the class of $M$ in $K_0(\mathfrak{M})$.

This categorical perspective overlooks the question of what combinatorial information about a matroid is encoded in its Tutte polynomial. Much of the inquiry into

1991 *Mathematics Subject Classification.* 05B35, 06A08.

Research performed while a General Member of MSRI during the 1996-1997 Program on Combinatorics, and supported by the Natural Sciences and Engineering Research Council of Canada under operating grant OGP0105392. Research at MSRI is supported in part by NSF grant DMS-9022140.

this subject has been motivated by this question, and various interpretations of coefficients or specializations of $T_M(x,y)$ have been given; see [4, 7, 10, 11, 12, 14, 16, 25, 26, 30, 31, 34] and elsewhere. In this paper we consider a specialization $Y_M(q,t)$ of $T_M(x,y)$ inspired by analogy with the Potts model from statistical mechanics. This connection, first made explicit by Kasteleyn and Fortuin [19], is implicit in equation (10) of Tutte [31], the last equation on page 331 of Oxley and Welsh [25], and elsewhere in the literature. Although $Y_M(q,t)$ is insensitive to the presence of loops, this is the only information lost in the specialization from $T_M(x,y)$ to $Y_M(q,t)$. Moreover, $Y_M(q,t)$ presents a great deal of information about $M$ in a very convenient form. For example, if $M$ is the graphic matroid of a graph $G=(V,E)$ with $c$ connected components, then $[q^{\#E}]t^c Y_M(q,t)$ is the chromatic polynomial of $G$, and for each natural number $i$, $[t^i]Y_M(1/2,t)$ is the probability that a random spanning edge-subgraph of $G$ has exactly $c+i$ connected components; Propositions 1.4 and 1.5 below give more general statements. (We use the notation $[z^a]F(z)$ to denote the coefficient of $z^a$ in the polynomial $F(z)$.)

The main purpose of this paper is to explain the meaning of the coefficients of $Y_M(1-p,t)$; Section 1 is intended to put this result in context. These coefficients are simply related to the ranks of the Whitney homology groups of the opposite partial orders of the independent set complexes of the duals of the truncations of $M$; the precise statements appear in Section 3. The proof is in two parts: in Section 2 we give a numerical formula for the coeffients of $Y_M(1-p,t)$ involving Möbius functions of certain partial orders associated with $M$, and in Section 3 we use established results of topological combinatorics to interpret these numbers as ranks of certain homology groups. In particular, we obtain as a corollary a new homological interpretation for the coefficients of the characteristic polynomial of a matroid.

## 1. The Potts model and the Tutte dichromate.

We begin by considering general finite graphs $G=(V,E)$, which may have both loops and multiple edges. The notations $n(G)$, $m(G)$, and $c(G)$ denote the number of vertices, edges, and connected components of $G$, respectively. The lattice of set partitions of $V$ is denoted by $\Pi_V$. For $e \in E$ and $\pi \in \Pi_V$, denote by $e \prec \pi$ the relation that there exists a block $B \in \pi$ such that both ends of $e$ are in $B$; as usual, $e \not\prec \pi$ denotes the negation of this relation. For $\pi \in \Pi_V$, let $\langle G:\pi \rangle := \#\{e \in E : e \not\prec \pi\}$. We define a polynomial $Z_G(q,t) \in \mathbb{Z}[q,t]$ by

$$(1.1) \qquad Z_G(q,t) := \sum_{\pi \in \Pi_V} q^{\langle G:\pi \rangle} t_{(\#\pi)}$$

in which $t_{(k)} := t(t-1)\cdots(t-k+1)$ is the $k$-th falling factorial polynomial. Clearly $Z_G(q,t)$ depends only upon the isomorphism class of $G$. As an example, Figure 1 depicts a graph $G$ and its associated polynomial $Z_G(q,t)$. When $t=N$ is a positive integer and $q=e^{J/kT}$, $Z_G(e^{J/kT},N)$ is (almost by definition) the partition function of the $N$-state Potts model on $G$ ($T$ is temperature, $k$ is the Boltzmann constant, and $J$ is a material parameter); see [3, 5, 17, 19, 20] for references and more information.

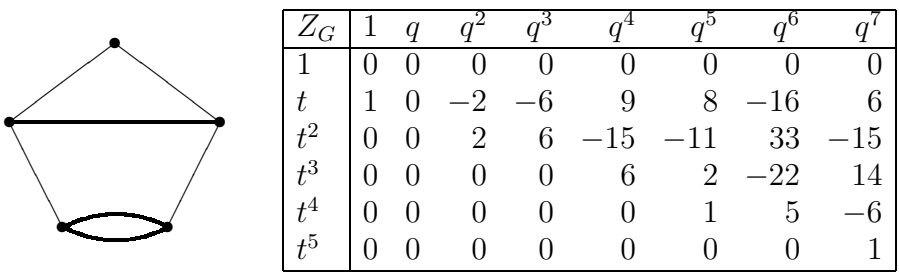

| $Z_G$ | 1 | $q$ | $q^2$ | $q^3$ | $q^4$ | $q^5$ | $q^6$ | $q^7$ |
|---|---|---|---|---|---|---|---|---|
| 1 | 0 | 0 | 0 | 0 | 0 | 0 | 0 | 0 |
| $t$ | 1 | 0 | $-2$ | $-6$ | 9 | 8 | $-16$ | 6 |
| $t^2$ | 0 | 0 | 2 | 6 | $-15$ | $-11$ | 33 | $-15$ |
| $t^3$ | 0 | 0 | 0 | 0 | 6 | 2 | $-22$ | 14 |
| $t^4$ | 0 | 0 | 0 | 0 | 0 | 1 | 5 | $-6$ |
| $t^5$ | 0 | 0 | 0 | 0 | 0 | 0 | 0 | 1 |

FIGURE 1. A graph $G$ and its polynomial $Z_G(q,t)$.

The polynomial $Z_G(q,t)$ can be calculated recursively as follows. For $G = (V, E)$ and $e \in E$, let $G \smallsetminus e$ and $G/e$ denote the graphs obtained from $G$ by deletion of $e$ and by contraction of $e$, respectively.

**Proposition 1.1.** *Let $G = (V, E)$ be a graph, and let $e \in E$. Then*
$$Z_G(q,t) = qZ_{G\smallsetminus e}(q,t) + (1-q)Z_{G/e}(q,t).$$
*If $G$ consists of $n$ vertices, any number of loops, and no nonloop edges, then $Z_G(q,t) = t^n$.*

*Proof.* The partitions $\pi \in \Pi_V$ such that $e \prec \pi$ are in a natural bijective correspondence with partitions in $\Pi_{V(G/e)}$. Thus we calculate that
$$\begin{aligned} Z_G(q,t) &= q \sum_{e \not\prec \pi \in \Pi_V} q^{\langle G\smallsetminus e:\pi\rangle} t_{(\#\pi)} + \sum_{e \prec \pi \in \Pi_V} q^{\langle G:\pi\rangle} t_{(\#\pi)} \\ &= qZ_{G\smallsetminus e}(q,t) + (1-q) \sum_{e \prec \pi \in \Pi_V} q^{\langle G\smallsetminus e:\pi\rangle} t_{(\#\pi)} \\ &= qZ_{G\smallsetminus e}(q,t) + (1-q)Z_{G/e}(q,t). \end{aligned}$$
The second assertion follows from (1.1) since $t^n = \sum_k S(n,k)t_{(k)}$, where $S(n,k)$ denotes a Stirling number of the second kind. □

**Proposition 1.2.** *If $G$ and $H$ are vertex-disjoint graphs then*
$$Z_{G \cup H}(q,t) = Z_G(q,t)Z_H(q,t).$$

*Proof.* Use Proposition 1.1 and induction on the number of nonloop edges of $G \cup H$. □

**Corollary 1.3.** *For any graph $G$, $t^{c(G)}$ divides $Z_G(q,t)$.*

*Proof.* From (1.1) it follows that $t$ divides $Z_C(q,t)$ for each connected component $C$ of $G$; the result follows from Proposition 1.2. □

The polynomial $Z_G(q,t)$ displays a great deal of information about $G$ in a very convenient form, by considering the rows and columns of Figure 1 (for example) as univariate polynomials. For $t \in \mathbb{N}$, a *t-colouring* of $G$ is any function $f : V \to \{1, 2, ..., t\}$. An edge $e \in E$ is *proper with respect to* $f$ if $e = \{u, v\}$ and $f(u) \neq f(v)$, so loops are never proper. Proposition 1.4 is equivalent to Theorem 3.2 of Brylawski [11]. (See also equations (10) and (12) of Tutte [31], and Theorem 6.3.26 of Brylawski and Oxley [12], in connection with Corollary 1.8 below.)

**Proposition 1.4.** *Let $G$ be a graph. For any $k, t \in \mathbb{N}$, $[q^k]Z_G(q,t)$ is the number of t-colourings of $G$ which have exactly $k$ proper edges.*

*Proof.* To each $t$-colouring $f$ of $G$ is associated the partition $\pi_f \in \Pi_V$ defined as follows: $v, w \in V$ are in the same block of $\pi_f$ if and only if $f(v) = f(w)$. For a given $\pi \in \Pi_V$, the number of $t$-colourings $f$ such that $\pi_f = \pi$ is $t_{(\#\pi)}$. The result now follows directly from (1.1). □

In particular, $[q^{m(G)}]Y_G(q,t)$ is the chromatic polynomial of $G$ [10, 11, 12, 18, 21, 26, 29, 30, 31, 35]. Given a graph $G$ and $0 \leq q \leq 1$, let $\mathcal{G}(q)$ denote a random spanning subgraph of $G$ obtained by deleting each edge independently with probability $q$. Proposition 1.5 (in conjunction with Corollary 1.8) is implicit in the last equation on page 331 of Oxley and Welsh [25] (and is derived and interpreted there more generally for matroids).

**Proposition 1.5.** *Let $G$ be a graph. For any $0 \leq q \leq 1$ and $k \in \mathbb{N}$, $[t^k]Z_G(q,t)$ is the probability that $\mathcal{G}(q)$ has exactly $k$ connected components.*

*Proof.* Fix any $k \in \mathbb{N}$ and $0 \leq q \leq 1$. We proceed by induction on the number of nonloop edges of $G$. The basis of induction follows from the second assertion of Proposition 1.1. For the induction step, let $e \in E$ be a nonloop, and denote by $\mathrm{rel}_k^G(q)$ the probability that $\mathcal{G}(q)$ has exactly $k$ connected components. By conditioning on the state of $e$ (deleted or included) we obtain

(1.2) $$\mathrm{rel}_k^G(q) = q\,\mathrm{rel}_k^{G \smallsetminus e}(q) + (1-q)\mathrm{rel}_k^{G/e}(q).$$

Comparison of Proposition 1.1 and (1.2) completes the induction step. □

In particular, $[t^1]Z_G(q,t)$ is the *reliability polynomial* of $G$, also studied under the title of percolation on graphs; see [2, 12, 13, 25, 32, 34].

For background information on matroids the reader should consult Oxley [24] or Welsh [33]. We view a matroid $M = (E, \rho)$ as a *ground-set* $E$ with a *rank function* $\rho : 2^E \to \mathbb{N}$ satisfying certain axioms; the *corank function* $\sigma : 2^E \to \mathbb{N}$ is defined by $\sigma(S) := \#S - \rho(S)$. The *rank* of $M$ is $\rho(E)$, and will also be denoted by $d(M)$. The simplicial complex of *independent sets* of $M$ is $\mathcal{I}(M) := \{S \subseteq E : \sigma(S) = 0\}$. To each graph $G = (V, E)$ is associated its *graphic* (or polygon, or cycle) matroid $M(G)$, which has ground-set $E$ and rank function $\rho(S) := n(G) - c((V, S))$; in this case $\mathcal{I}(M)$ consists of all edge-sets of spanning forests of $G$. For a matroid $M = (E, \rho)$ and $e \in E$ we denote by $M \smallsetminus e$ and $M/e$ the matroids obtained from $M$ by deletion of $e$ and by contraction of $e$, respectively. For $S \subseteq E$ and $e \in S$, let $S/e$ denote $S \smallsetminus e$ regarded as a subset of $E(M/e)$. We regard $S \smallsetminus e$ either as a subset of $E(M)$ or of $E(M \smallsetminus e)$, as required. The elements of $E$ are of three types, called *loops*, *coloops*, and *links*, and the following relations on ranks and coranks hold.

(1.3)

| type of $e$ | $\rho(E \smallsetminus e)$ | $\sigma(E \smallsetminus e)$ | $\rho(E/e)$ | $\sigma(E/e)$ |
|---|---|---|---|---|
| loop | $\rho(E)$ | $\sigma(E) - 1$ | $\rho(E)$ | $\sigma(E) - 1$ |
| coloop | $\rho(E) - 1$ | $\sigma(E)$ | $\rho(E) - 1$ | $\sigma(E)$ |
| link | $\rho(E)$ | $\sigma(E) - 1$ | $\rho(E) - 1$ | $\sigma(E)$ |

The *Tutte polynomial* of $M = (E, \rho)$, as defined by Crapo [14], is

$$(1.4) \qquad T_M(x, y) := \sum_{S \subseteq E} (x - 1)^{\rho(E) - \rho(S)} (y - 1)^{\sigma(S)}.$$

Proposition 1.1 shows that for any finite graph $G$, $Z_G(q, t) = f(G; t, 1 - q, q)$ in which $f(G; t, x, y)$ is a polynomial invariant of graphs defined by Negami [21]. Oxley [23] shows that $t^{-c(G)} f(G; t, x, y)$ is determined by the Tutte polynomial of $M(G)$. This suffices to establish Corollary 1.8 below, but for completeness we sketch another proof.

For a minor-closed class $\mathfrak{C}$ of matroids, a *Tutte-Grothendieck invariant* on $\mathfrak{C}$ is an assignment $M \mapsto \Psi_M$ which associates to each matroid $M$ in $\mathfrak{C}$ an element $\Psi_M$ of some fixed commutative ring $R$ (which depends only on $\Psi$), and satisfies the following properties:
($\alpha$) If $M \simeq N$ then $\Psi_M = \Psi_N$.
($\beta$) If $e \in E(M)$ is a link then $\Psi_M = \Psi_{M \smallsetminus e} + \Psi_{M/e}$.
($\gamma$) If $M$ and $N$ have disjoint ground-sets then $\Psi_{M \oplus N} = \Psi_M \Psi_N$.
For a graph $G$, let $r(G) := n(G) - c(G)$ and let $s(G) := m(G) - n(G) + c(G)$. Define a rational function $\widetilde{Z}_G(q, t) \in \mathbb{Z}(q, t)$ by

$$(1.5) \qquad \widetilde{Z}_G(q, t) := \frac{Z_G(q, t)}{t^{c(G)} (1 - q)^{r(G)} q^{s(G)}}.$$

**Proposition 1.6.** *The assignment $M(G) \mapsto \widetilde{Z}_G(q, t)$ is a (well-defined) Tutte-Grothendieck invariant on the class of graphic matroids.*

*Sketch of proof.* Well-definedness is the main issue. Whitney [36] (see also Theorem 6.3.1 of [22]) proves that two graphs $G$ and $H$ are such that $M(G) \simeq M(H)$ if and only if $G$ and $H$ are "2-isomorphic" (see the above references for the definition). If $H$ is obtained from $G$ by splitting at a cut-vertex or by twisting at a 2-vertex-cut then Proposition 1.1 and induction on $m(G)$ can be used to show that $Z_G(q, t) = Z_H(q, t)$. This suffices to show that if $G$ and $H$ are 2-isomorphic then $Z_G(q, t) = Z_H(q, t)$, and it follows that $M(G) \mapsto \widetilde{Z}_G(q, t)$ is well-defined. Property ($\alpha$) is a consequence of well-definedness, ($\beta$) follows from Proposition 1.1, and ($\gamma$) follows from Proposition 1.2. □

As mentioned in the introduction, $T_M(x,y)$ is the *universal Tutte-Grothendieck invariant* (on the class of all matroids) in the following sense. Let $L$ denote a one-element matroid of rank zero, and let $L^*$ denote a one-element matroid of rank one; notice that $T_{L^*}(x,y) = x$ and $T_L(x,y) = y$. Proposition 1.7 is developed further in [10, 12, 25, 37].

**Proposition 1.7.** *If $\mathfrak{C}$ is a minor-closed class of matroids and $M \mapsto \Psi_M$ is a Tutte-Grothendieck invariant on $\mathfrak{C}$, then for all $M$ in $\mathfrak{C}$, $\Psi_M = T_M(\Psi_{L^*}, \Psi_L)$.*

*Sketch of proof.* It follows from ($\alpha$), ($\beta$), and ($\gamma$), by induction on $\#E(M)$, that the assignment $M \mapsto \Psi_M$ is determined by its values $\Psi_{L^*}$ and $\Psi_L$. Since $x$ and $y$ are independent indeterminates it follows that $\Psi_M = T_M(\Psi_{L^*}, \Psi_L)$ once it has been shown that $M \mapsto T_M(x,y)$ is itself a Tutte-Grothendieck invariant on $\mathfrak{M}$. Verification of this fact follows from the definition (1.4) as in Lemma 6.2.1 of [12], for example. □

See [19, 25, 31] for more information related to Corollary 1.8.

**Corollary 1.8.** *Let $G$ be a graph, with graphic matroid $M$. Then*
$$Z_G(q,t) = t^{c(G)}(1-q)^{r(G)}q^{s(G)} T_M\left(\frac{qt+1-q}{1-q}, \frac{1}{q}\right).$$

*Proof.* One checks that $\widetilde{Z}_{L^*}(q,t) = (qt+1-q)/(1-q)$ and $\widetilde{Z}_L = 1/q$, and the result follows from Propositions 1.6 and 1.7. □

Corollary 1.8 allows us to generalize from the class of graphic matroids to the class of all matroids. However, $c(G)$ can not be determined from $M(G)$ alone, and so we define polynomials

(1.6) $$Y_G(q,t) := t^{-c(G)}Z_G(q,t)$$

for every graph $G$, and

(1.7) $$Y_M(q,t) := (1-q)^{\rho(E)}q^{\sigma(E)} T_M\left(\frac{qt+1-q}{1-q}, \frac{1}{q}\right)$$

for every matroid $M = (E, \rho)$. Corollary 1.8 shows that these definitions agree for graphic matroids. From ($\alpha$), ($\beta$), and ($\gamma$) we have the following.
($\alpha'$) If $M \simeq N$ then $Y_M(q,t) = Y_N(q,t)$.
($\beta'$) If $e \in E(M)$ is a link then $Y_M(q,t) = qY_{M \smallsetminus e}(q,t) + (1-q)Y_{M/e}(q,t)$.
($\gamma'$) If $M$ and $N$ have disjoint ground-sets then $Y_{M \oplus N}(q,t) = Y_M(q,t)Y_N(q,t)$.
If $e$ is a loop of $M$ then $M \simeq (M \smallsetminus e) \oplus L$; from ($\alpha$), ($\beta$), and (1.7) we deduce that
($\beta''$) If $e \in E(M)$ is a loop then $Y_M(q,t) = Y_{M \smallsetminus e}(q,t)$.
If $e$ is a coloop of $M$ then $M \simeq (M/e) \oplus L^*$; from ($\alpha$), ($\beta$), and (1.7) we deduce that
($\beta^*$) If $e \in E(M)$ is a coloop then $Y_M(q,t) = (qt+1-q)Y_{M/e}(q,t)$.

| $Y_M$ | 1 | $p$ | $p^2$ | $p^3$ | $p^4$ | $p^5$ | $p^6$ | $p^7$ |
|---|---|---|---|---|---|---|---|---|
| 1 | 0 | 0 | 0 | 0 | 19 | $-38$ | 26 | $-6$ |
| $t$ | 0 | 0 | 0 | 29 | $-100$ | 128 | $-72$ | 15 |
| $t^2$ | 0 | 0 | 20 | $-94$ | 176 | $-164$ | 76 | $-14$ |
| $t^3$ | 0 | 7 | $-41$ | 100 | $-130$ | 95 | $-37$ | 6 |
| $t^4$ | 1 | $-7$ | 21 | $-35$ | 35 | $-21$ | 7 | $-1$ |

FIGURE 2. The polynomial $Y_M(1-p,t)$ for $M(G)$ from Figure 1.

| $\widehat{Y}_M$ | 1 | $p$ | $p^2$ | $p^3$ | $p^4$ | $p^5$ | $p^6$ | $p^7$ |
|---|---|---|---|---|---|---|---|---|
| 1 | 0 | 0 | 0 | 0 | 19 | 38 | 26 | 6 |
| $t$ | 0 | 0 | 0 | 29 | 81 | 90 | 46 | 9 |
| $t^2$ | 0 | 0 | 20 | 65 | 95 | 74 | 30 | 5 |
| $t^3$ | 0 | 7 | 21 | 35 | 35 | 21 | 7 | 1 |

FIGURE 3. The polynomial $\widehat{Y}_M(p,t)$ for $M(G)$ from Figure 1.

Let $U_m^0$ be a matroid of rank zero on $m$ elements; then $T_{U_m^0}(x,y) = y^m$, but $Y_{U_m^0}(q,t) = 1$ is independent of $m$. Generally, $Y_M(q,t)$ is insensitive to the presence of loops, but this is the only information lost in specializing from $T_M(x,y)$ to $Y_M(q,t)$. To see this, invert (1.7) to yield

$$(1.8) \qquad T_M(x,y) = \frac{y^{\#E}}{(y-1)^{\rho(E)}} Y_M\left(\frac{1}{y}, (x-1)(y-1)\right).$$

The degree of $Y_M(q,t)$ in the variable $t$ is $\rho(E)$, and the degree of $Y_M(q,t)$ in the variable $q$ is the number of nonloop elements of $M$. Thus, given $Y_M(q,t)$ and the number of loops of $M$, $T_M(x,y)$ can be recovered.

## 2. THE COEFFICIENTS OF $Y_M(1-p,t)$.

We now make the substitution $p := 1-q$ and interpret the coefficients of $Y_M(1-p,t)$ combinatorially. Figure 2 shows this form of the polynomial for the graph from Figure 1. One notices in this example a simple sign-alternation pattern, and the fact that evaluation at $t=1$ yields $Y_M(1-p,1) = 1$. This suggests examination of the rational function $\widehat{Y}_M(p,t) \in \mathbb{Z}(p,t)$ defined by

$$(2.1) \qquad \widehat{Y}_M(p,t) := \frac{(-1)^{d(M)} Y_M(1+p,-t) - t^{d(M)}}{1+t}.$$

Figure 3 shows $\widehat{Y}_M(p,t)$ for the example from Figures 1 and 2. For $e \in E(M)$, define

$b(M, e) := \rho(E) - \rho(E \smallsetminus e)$, so $b(M, e) = 1$ if $e$ is a coloop of $M$ and $b(M, e) = 0$ otherwise.

**Lemma 2.1.** *Let $M = (E, \rho)$ be a matroid and let $e \in E$ be a link or a coloop. Then*
$$\widehat{Y}_M(p, t) = pt^{d(M)-1} + (1+p)t^{b(M,e)}\widehat{Y}_{M \smallsetminus e}(p, t) + p\widehat{Y}_{M/e}(p, t).$$
*Consequently, for every matroid $M$, $\widehat{Y}_M(p, t) \in \mathbb{N}[p, t]$.*

*Proof.* Let $b := b(M, e)$ and $d := d(M)$. From $(\beta')$ and $(\beta^*)$ and the fact that if $e \in E$ is a coloop then $M \smallsetminus e \simeq M/e$ we see that
$$Y_M(q, t) = qt^b Y_{M \smallsetminus e}(q, t) + (1 - q)Y_{M/e}(q, t).$$
Therefore,
$$\begin{aligned}
\widehat{Y}_M(p, t) &= \frac{(-1)^d(1+p)(-t)^b Y_{M \smallsetminus e}(1+p, -t)}{1+t} \\
&\quad + \frac{(-1)^d(-p)Y_{M/e}(1+p, -t) - t^d}{1+t} \\
&= (1+p)t^b \widehat{Y}_{M \smallsetminus e}(p, t) + p\widehat{Y}_{M/e}(p, t) + pt^{d-1},
\end{aligned}$$
as claimed. The second claim follows by induction on the number of nonloop elements of $E$; for the basis of induction $d = 0$ and $Y_M(q, t) = 1$, so that $\widehat{Y}_M(p, t) = 0$. □

Let $M = (E, \rho)$ be a matroid of rank $d \geq 1$. For each $0 \leq i \leq d - 1$ let $\mathcal{S}_i^M$ be the set of all $S \subseteq E$ such that $\rho(S) \geq d - i$. Notice that $\mathcal{S}_0^M \subset \mathcal{S}_1^M \subset \cdots \subset \mathcal{S}_{d-1}^M$, and that if $S \subseteq T \subseteq E$ and $S \in \mathcal{S}_i^M$ then $T \in \mathcal{S}_i^M$. For example, if $M = M(G)$ for a graph $G$ then $d = n(G) - c(G)$ and for $0 \leq i \leq d - 1$, $\mathcal{S}_i^M$ consists of the edge-sets of all spanning subgraphs of $G$ which have at most $c(G) + i$ components. Order $\mathcal{S}_i^M$ partially by inclusion, and adjoin a new unique minimal element $\widehat{0}$ to $\mathcal{S}_i^M$ to obtain $\mathcal{L}_i^M := \{\widehat{0}\} \oplus \mathcal{S}_i^M$, which is in fact a graded lattice. The height of $S \in \mathcal{S}_i^M$ in $\mathcal{L}_i^M$ is $\#S - d + 1 + i$. Let $\mu_i^M(\cdot, \cdot)$ denote the Möbius function of $\mathcal{L}_i^M$ (see Rota [26]); we recall that $\mu_i^M(\widehat{0}, \widehat{0}) = 1$ and that for $S \in \mathcal{S}_i^M$,

(2.2) $$\mu_i^M(\widehat{0}, S) = -\sum_{\widehat{0} \leq T < S} \mu_i^M(\widehat{0}, T).$$

If $S \subseteq E$ but $S \notin \mathcal{S}_i^M$ then by convention we put $\mu_i^M(\widehat{0}, S) := 0$. Throughout this paper the first argument of a Möbius function is always $\widehat{0}$, and we shall henceforth omit it from the notation, writing only $\mu_i^M(S)$.

**Lemma 2.2.** *Let $M = (E, \rho)$ be a matroid of rank $d \geq 1$ and let $e \in E$. Then for all $0 \leq i \leq d - 1$ and all $S \in \mathcal{S}_i^M$ with $e \in S$:*
(a) *if $e$ is a loop then $\mu_i^M(S) = 0$;*

(b) *if $e$ is a coloop then $\mu_i^M(S) = -\mu_{i-1}^{M\smallsetminus e}(S\smallsetminus e) + \mu_i^{M/e}(S/e)$;*
(c) *if $e$ is a link then $\mu_i^M(S) = -\mu_i^{M\smallsetminus e}(S\smallsetminus e) + \mu_i^{M/e}(S/e)$.*

*Proof.* We prove the claims by induction on $\#S$. For the basis of induction notice that if $S \in \mathcal{S}_i^M$ then $S$ contains an independent set of rank $d-i$, and so $\#S \geq d-i$. If $e$ is a loop then this bound can not be attained, and the basis of induction is that $S\smallsetminus e \in \mathcal{I}(M)$ has rank $d-i$. Thus, the interval from $\widehat{0}$ to $S$ in $\mathcal{L}_i^M$ consists of the three elements $\{\widehat{0} < S\smallsetminus e < S\}$, and so $\mu_i^M(S) = 0$, as required. If $e$ is a coloop and $S \in \mathcal{I}(M)$ has rank $d-i$ then $S\smallsetminus e \in \mathcal{I}(M\smallsetminus e)$ has rank $(d-1)-(i-1)-1$, so $S\smallsetminus e \notin \mathcal{S}_{i-1}^{M\smallsetminus e}$, so $\mu_{i-1}^{M\smallsetminus e}(S\smallsetminus e) = 0$. If $e$ is a link and $S \in \mathcal{I}(M)$ has rank $d-i$ then $S\smallsetminus e \in \mathcal{I}(M\smallsetminus e)$ has rank $d-i-1$, so $S\smallsetminus e \notin \mathcal{S}_i^{M\smallsetminus e}$, so $\mu_i^{M\smallsetminus e}(S\smallsetminus e) = 0$. In both these cases, $S$ and $S/e$ cover $\widehat{0}$ in $\mathcal{L}_i^M$ and $\mathcal{L}_i^{M/e}$, respectively, so that $\mu_i^M(S) = \mu_i^{M/e}(S/e) = -1$, as required. This completes the basis of induction.

For the induction step apply (2.2) and separate the case $e \notin T$ from the case $e \in T$ to obtain

$$(2.3) \qquad \mu_i^M(S) = -\sum_{T \leq S\smallsetminus e} \mu_i^M(T) - \sum_{e \in T < S} \mu_i^M(T).$$

The first sum on the right side of (2.3) is zero, since it is the sum of $\mu_i^M(T)$ over the interval $[\widehat{0}, S\smallsetminus e]$ of $\mathcal{L}_i^M$. We may apply the induction hypothesis to each term in the second sum. If $e$ is a loop then each term of this sum is zero, proving part (a). If $e$ is a coloop or a link then let $b := b(M, e)$. By induction and (2.3) we have

$$\mu_i^M(S) = -\sum_{e \in T < S} \left( -\mu_{i-b}^{M\smallsetminus e}(T\smallsetminus e) + \mu_i^{M/e}(T/e) \right)$$
$$= \sum_Q \mu_{i-b}^{M\smallsetminus e}(Q) - \sum_U \mu_i^{M/e}(U),$$

in which $Q$ ranges over subsets of $E\smallsetminus e$ such that $Q \cup \{e\} \in \mathcal{L}_i^M$ and $Q < S\smallsetminus e$, and $U$ ranges over subsets of $E/e$ such that $U \cup \{e\} \in \mathcal{L}_i^M$ and $U < S/e$. If $i = d-1$ then $Q = \varnothing$ and $U = \varnothing$ arise, and can be identified with the terms for $\widehat{0}$ in $\mathcal{L}_{i-b}^{M\smallsetminus e}$ and in $\mathcal{L}_i^{M/e}$, respectively. If $i \leq d-2$ then these terms do not arise, but since $\mu_{i-b}^{M\smallsetminus e}(\widehat{0}) - \mu_i^{M/e}(\widehat{0}) = 0$ we may include them without changing the sum. We obtain

$$\mu_i^M(S) = \sum_{\widehat{0} \leq Q < S\smallsetminus e} \mu_{i-b}^{M\smallsetminus e}(Q) - \sum_{\widehat{0} \leq U < S/e} \mu_i^{M/e}(U)$$
$$= -\mu_{i-b}^{M\smallsetminus e}(S\smallsetminus e) + \mu_i^{M/e}(S/e),$$

by (2.2). This completes the induction step, and the proof. $\square$

**Corollary 2.3.** *Let $M = (E, \rho)$ be a matroid of rank $d \geq 1$, let $0 \leq i \leq d-1$, and let $S \in \mathcal{S}_i^M$. Then $(-1)^{\#S-d+1+i}\mu_i^M(S) \geq 0$. Moreover, if $S$ contains no loops of $M$ then $(-1)^{\#S-d+1+i}\mu_i^M(S) > 0$.*

*Proof.* Both claims follows from Lemma 2.2 by induction on $\#S$. □

Let $M = (E, \rho)$ be a matroid of rank $d \geq 1$, with $\#E = m$ and with $\ell$ loops. For $0 \leq i \leq d-1$ define polynomials $W_i^M(p) \in \mathbb{Z}[p]$ by

$$W_i^M(p) := \sum_{S \in \mathcal{S}_i^M} \mu_i^M(S)(-p)^{\#S-d+1+i}. \tag{2.4}$$

The coefficients $\omega_j^{i,M}$ defined by $W_i^M(p) = \sum_{j=1}^{m-d+1+i} \omega_j^{i,M} p^j$ are the *proper unsigned Whitney numbers of the first kind* of $\mathcal{L}_i^M$; that is,

$$\omega_j^{i,M} := (-1)^j \sum_{S \in \mathcal{S}_i^M:\ \#S-d+1+i=j} \mu_i^M(S) \tag{2.5}$$

for each $0 \leq i \leq d-1$ and $1 \leq j \leq m-d+1+i$. By convention we put $\omega_j^{i,M} := 0$ in all remaining cases, and from Corollary 2.3 it follows that $\omega_j^{i,M} \geq 0$ for all $i$ and $j$. By Lemma 2.2, $\mu_i^M(S) = 0$ for all $0 \leq i \leq d-1$ whenever $S \subseteq E$ contains a loop of $M$. From this, Corollary 2.3, and (2.4), it follows that $\deg W_i^M(p) = m - \ell - d + 1 + i$ for all $0 \leq i \leq d-1$.

**Lemma 2.4.** *Let $M = (E, \rho)$ be a matroid of rank $d \geq 1$ with $\#E = m$ and with $\ell$ loops. Then $W_{d-1}^M(p) = (1+p)^{m-\ell} - 1$.*

*Proof.* By (2.4) and Lemma 2.2 we have

$$W_{d-1}^M(p) = \sum_{S \in \mathcal{T}} \mu_{d-1}^M(S)(-p)^{\#S},$$

in which $\mathcal{T}$ is the set of all $S \in \mathcal{S}_{d-1}^M$ which do not contain any loops of $M$. One checks that $\mathcal{T}$ is a downward-closed subset of $\mathcal{S}_{d-1}^M$ and is the set of all nonempty subsets of the set of nonloop elements of $M$. Thus $\{\widehat{0}\} \oplus \mathcal{T}$ is a Boolean lattice with $m - \ell$ atoms, and the result follows since $\mu_{d-1}^M(S) = (-1)^{\#S}$ for all $S \in \mathcal{T}$. □

**Lemma 2.5.** *Let $M = (E, \rho)$ be a matroid of rank $d \geq 2$, and let $0 \leq i \leq d-2$. If $e \in E$ is a link or a coloop then*

$$W_i^M(p) = (1+p)W_{i-b(M,e)}^{M \smallsetminus e}(p) + W_i^{M/e}(p).$$

*Proof.* Let $b := b(M, e)$. Notice that $\{S \in \mathcal{S}_i^M : e \notin S\} = \{S \subseteq E(M \smallsetminus e) : \rho(S) \geq (d-b) - (i-b)\} = \mathcal{S}_{i-b}^{M \smallsetminus e}$. Also, for all $S$ in this set $\mu_i^M(S) = \mu_{i-b}^{M \smallsetminus e}(S)$ since the intervals $[\widehat{0}, S]$ in $\mathcal{L}_i^M$ and in $\mathcal{L}_{i-b}^{M \smallsetminus e}$ are isomorphic. Since $i \leq d-2$ the singleton $\{e\}$ is not an element of $\mathcal{S}_i^M$. Hence, $\{S \smallsetminus e : e \in S \in \mathcal{S}_i^M\} = \{Q \subseteq E(M \smallsetminus e) : \rho(Q) \geq$

$(d-b) - (i-b+1)\} = \mathcal{S}_{i-b+1}^{M \smallsetminus e} \supset \mathcal{S}_{i-b}^{M \smallsetminus e}$. Also, $\mu_{i-b}^{M \smallsetminus e}(Q) = 0$ if $Q \in \mathcal{S}_{i-b+1}^{M \smallsetminus e} \smallsetminus \mathcal{S}_{i-b}^{M \smallsetminus e}$. Finally, $\{S/e : e \in S \in \mathcal{S}_i^M\} = \{U \subseteq E(M/e) : \rho(U) \geq (d-1) - i\} = \mathcal{S}_i^{M/e}$. Thus, substituting Lemma 2.2 into (2.4) we obtain

$$
\begin{aligned}
W_i^M(p) &= \sum_{e \notin S \in \mathcal{S}_i^M} \mu_i^M(S)(-p)^{\#S - d + 1 + i} \\
&\quad + \sum_{e \in S \in \mathcal{S}_i^M} \left( -\mu_{i-b}^{M \smallsetminus e}(S \smallsetminus e) + \mu_i^{M/e}(S/e) \right) (-p)^{\#S - d + 1 + i} \\
&= \sum_{S \in \mathcal{S}_{i-b}^{M \smallsetminus e}} \mu_{i-b}^{M \smallsetminus e}(S)(-p)^{\#S - (d-b) + 1 + (i-b)} \\
&\quad - \sum_{Q \in \mathcal{S}_{i-b+1}^{M \smallsetminus e}} \mu_{i-b}^{M \smallsetminus e}(Q)(-p)^{(\#Q + 1) - (d-b) + 1 + (i-b)} \\
&\quad + \sum_{U \in \mathcal{S}_i^{M/e}} \mu_i^{M/e}(U)(-p)^{\#U - (d-1) + 1 + i} \\
&= W_{i-b}^{M \smallsetminus e}(p) + p W_{i-b}^{M \smallsetminus e}(p) + W_i^{M/e}(p),
\end{aligned}
$$

as required. $\square$

**Theorem 2.6.** *Let $M = (E, \rho)$ be a matroid of rank $d \geq 0$. Then*

$$\widehat{Y}_M(p, t) = \sum_{i=0}^{d-1} W_i^M(p) p^{d-1-i} t^i.$$

*Proof.* We prove, by induction on the number of nonloops of $M$, that for all $0 \leq i \leq d-1$, $[t^i]\widehat{Y}_M(p, t) = W_i^M(p) p^{d-1-i}$. If $M$ has no nonloops then $d(M) = 0$, $\widehat{Y}_M(p, t) = 0$, and the summation on the right side is empty; this establishes the basis of induction. For the induction step let $e$ be a coloop or link of $M$, let $b := b(M, e)$, and let $M$ have $\ell$ loops and $\#E(M) = m$. The induction hypothesis is that $[t^i]\widehat{Y}_{M \smallsetminus e}(p, t) = W_i^{M \smallsetminus e}(p) p^{(d-b)-1-i}$ for all $0 \leq i \leq d-b-1$ and that $[t^i]\widehat{Y}_{M/e}(p, t) = W_i^{M/e}(p) p^{(d-1)-1-i}$ for all $0 \leq i \leq d-2$. For the case $i = d-1$ Lemmas 2.1, 2.4, and the induction hypothesis imply that

$$
\begin{aligned}
[t^{d-1}]\widehat{Y}_M(p, t) &= p + (1+p)[t^{(d-b)-1}]\widehat{Y}_{M \smallsetminus e}(p, t) + p[t^{d-1}]\widehat{Y}_{M/e}(p, t) \\
&= p + (1+p)\left((1+p)^{m-1-\ell} - 1\right) + p \cdot 0 \\
&= (1+p)^{m-\ell} - 1 = W_{d-1}^M(p) p^{d-1-(d-1)},
\end{aligned}
$$

as required. For the case $0 \leq i \leq d-2$ Lemmas 2.1, 2.5, and the induction hypothesis imply that

$$\begin{aligned}
[t^i]\widehat{Y}_M(p,t) &= (1+p)[t^{i-b}]\widehat{Y}_{M\smallsetminus e}(p,t) + p[t^i]\widehat{Y}_{M/e}(p,t) \\
&= (1+p)W_{i-b}^{M\smallsetminus e}(p)p^{(d-b)-1-(i-b)} + pW_i^{M/e}(p)p^{d-2-i} \\
&= \left((1+p)W_{i-b}^{M\smallsetminus e}(p) + W_i^{M/e}(p)\right)p^{d-1-i} \\
&= W_i^M(p)p^{d-1-i},
\end{aligned}$$

completing the induction step, and the proof. □

Inverting (2.1), we see that

(2.6) $$Y_M(1-p,t) = t^{d(M)} + (1-t)(-1)^{d(M)}\widehat{Y}_M(-p,-t).$$

Corollary 2.7 is an immediate consequence of Theorem 2.6 and (2.6).

**Corollary 2.7.** *Let $M=(E,\rho)$ be a matroid of rank $d \geq 1$ with $\#E = m$ and with $\ell$ loops. For all $0 \leq i \leq d$ and $0 \leq k \leq m-\ell$,*

$$[p^k t^i]Y_M(1-p,t) = \delta_{k0}\delta_{id} + (-1)^{k+i-d}\left(\omega_{k+i-d+1}^{i,M} + \omega_{k+i-d}^{i-1,M}\right),$$

*in which $\delta_{xy}$ is the Kronecker delta function, and $[p^k t^i]Y_M(1-p,t) = 0$ for all other values of $k$ and $i$.*

**Corollary 2.8.** *Let $M=(E,\rho)$ be a matroid of rank $d \geq 1$. For $0 \leq k \leq d$, the number of $k$-element sets in $\mathcal{I}(M)$ is $[p^k t^{d-k}]Y_M(1-p,t)$.*

*Proof.* For each $1 \leq i \leq d-1$, the minimal elements of $\mathcal{S}_i^M$ are the sets $S \in \mathcal{I}(M)$ such that $\#S = d-i$. Since each such set has $\mu_i^M(S) = -1$, (2.5) shows that the number of $(d-i)$-element sets in $\mathcal{I}(M)$ is $\omega_1^{i,M}$ for each $1 \leq i \leq d-1$. From Corollary 2.7 we obtain

$$[p^k t^{d-k}]Y_M(1-p,t) = \delta_{k0} + \omega_1^{d-k,M} + \omega_0^{d-k-1,M}$$

for all $0 \leq k \leq d$, and the result follows. □

For a matroid $M$, let $\mathcal{K}(M)$ denote the geometric lattice of closed sets (or "flats") of $M$. The *characteristic polynomial* of $M$ is

$$\chi_M(t) := \sum_{S \in \mathcal{K}(M)} \mu_{\mathcal{K}}(S) t^{d(M)-\rho(S)}.$$

See Rota [26] for more information. In particular, $\chi_L(t) = 1$ and $\chi_{L^*}(t) = t-1$, if $e \in E(M)$ is a link then $\chi_M(t) = \chi_{M\smallsetminus e}(t) - \chi_{M/e}(t)$, and if $M$ and $N$ have disjoint ground-sets then $\chi_{M \oplus N}(t) = \chi_M(t)\chi_N(t)$. From this it follows by induction on $m := \#E(M)$ that if $M$ has $\ell$ loops and rank $d$, then

$$\chi_M(t) = [q^{m-\ell}]Y_M(q,t) = (-1)^{m-\ell+d}[p^{m-\ell}]Y_M(1-p,t).$$

Corollary 2.9 follows directly from (2.6).

**Corollary 2.9.** *Let $M = (E, \rho)$ be a matroid of rank $d \geq 1$ with $\#E = m$ and with no loops. The characteristic polynomial of $M$ is*

$$\chi_M(t) = (-1)^d (1-t) \sum_{i=0}^{d-1} \omega_{m-d+1+i}^{i,M} (-t)^i.$$

## 3. Whitney homology of upside-down matroids.

For this section we assume familiarity with the basic concepts of algebraic topology applied to combinatorics, as surveyed by Björner [8]. The *order complex* of a bounded partial order $P$ with $\widehat{0} \neq \widehat{1}$ is the set $\Delta(P)$ of all totally-ordered subsets of $P \smallsetminus \{\widehat{0}, \widehat{1}\}$, partially ordered by inclusion; $\Delta(P)$ is an abstract simplicial complex. The (integral reduced simplicial) homology groups of a $g$-dimensional simplicial complex $\Delta$ with coefficients in $\mathbb{Z}$ are denoted by $\widetilde{H}_*(\Delta; \mathbb{Z}) = \bigoplus_{j=-1}^{g} \widetilde{H}_j(\Delta; \mathbb{Z})$. The *Betti numbers* of $\Delta$ are $\beta_j(\Delta) := \operatorname{rank} \widetilde{H}_j(\Delta; \mathbb{Z})$ for $-1 \leq j \leq g$. (For $\Delta \neq \varnothing$, $\beta_{-1}(\Delta) = 0$; also $\beta_{-1}(\varnothing) = 1$ and $\beta_j(\varnothing) = 0$ for $j > 0$.) Philip Hall's Theorem and the Euler-Poincaré Formula relate the Möbius funtion of $P$ to the (reduced) Euler characteristic of $\Delta(P)$:

$$\mu_P(\widehat{1}) = \widetilde{\chi}(\Delta(P)) = -\beta_{-1} + \beta_0 - \beta_1 + \cdots + (-1)^g \beta_g.$$

For a partial order $P$ with $\widehat{0}$ the *Whitney homology* $WH_*(P; \mathbb{Z})$ of $P$ was defined by Baclawski [1]. Theorem 5.1 of Björner [6] shows that $WH_*(P; \mathbb{Z}) = \bigoplus_{j=0}^{g} WH_j(P; \mathbb{Z})$ in which

$$WH_j(P; \mathbb{Z}) = \bigoplus_{\widehat{0} < \sigma \in P} \widetilde{H}_{j-1}(\Delta([\widehat{0}, \sigma]); \mathbb{Z})$$

for each $0 \leq j \leq g$.

Interest in Whitney homology of geometric lattices arose as it provides an interpretation for the coefficients of the characteristic polynomial of a matroid; see [1, 6, 7]. Also, Theorem 9 of Stanley [28] implicitly uses Whitney homology of the opposite partial order of a geometric lattice to describe the Betti numbers of the minimal free resolution of the face-ring of a matroid; see also Eagon and Reiner [15]. Here we interpret the coefficients of $\widehat{Y}_M(p, t)$ as the ranks of the Whitney homology groups of the lattices $\mathcal{L}_i^M$, which are in general atomic but not semimodular.

Given a matroid $M = (E, \rho)$ of rank $d \geq 1$, an index $0 \leq i \leq d-1$, and a set $S \in \mathcal{S}_i^M$, let $\Gamma_i(S) := \{S \smallsetminus T : T \in \mathcal{S}_i^M \text{ and } T \subseteq S\}$; this simplicial complex depends only upon $i$ and the restriction submatroid $M|_S := M \smallsetminus (E \smallsetminus S)$ of $M$. In fact, denoting by $N_i$ the matroid dual to the $i$-th truncation of $M$, we have $\Gamma_i(S) = \mathfrak{I}(N_i/(E \smallsetminus S))$. On the other hand, the interval $[\widehat{0}, S]$ in $\mathcal{L}_i^M$ is isomorphic to the partial order opposite to $\Gamma_i(S) \oplus \{\widehat{1}\}$, and hence $\Delta([\widehat{0}, S]_{\mathcal{L}_i^M})$ is the barycentric subdivision of $\Gamma_i(S)$. This gives us enough information to understand the Whitney homology of the lattices $\mathcal{L}_i^M$ for $0 \leq i \leq d-1$.

**Theorem 3.1.** *Let $M = (E, \rho)$ be a matroid of rank $d \geq 1$ with $\#E = m$ and with no loops. For all $0 \leq i \leq d-1$ and $1 \leq j \leq m - d + 1 + i$,*
$$[p^{d-1-i+j} t^i] \widehat{Y}_M(p,t) = \omega_j^{i,M} = \operatorname{rank} WH_{j-1}(\mathcal{L}_i^M; \mathbb{Z}).$$

*Proof.* For each $0 \leq i \leq d-1$ let $N_i$ be the matroid dual to the $i$-th truncation of $M$, as above. For each $S \in \mathcal{S}_i^M$, $N_i/(E \smallsetminus S)$ is a matroid, and so $\Gamma_i(S) = \mathcal{I}(N_i/(E \smallsetminus S))$ is a shellable simplicial complex (Theorem 7.3.3 of [7]). Hence, only the highest-dimensional homology group of $\Gamma_i(S)$ may be nonzero, and it is torsion-free (Theorem 7.8.1 of [7]); this highest dimension is $\#S - d + i - 1$. Since $\Delta([\widehat{0}, S]_{\mathcal{L}_i^M})$ is the barycentric subdivision of $\Gamma_i(S)$ we have
$$\widetilde{H}_*(\Delta([\widehat{0}, S]_{\mathcal{L}_i^M}); \mathbb{Z}) = \widetilde{H}_*(\Gamma_i(S); \mathbb{Z}).$$

It follows that
$$\mu_i^M(S) = \widetilde{\chi}(\Delta([\widehat{0}, S]_{\mathcal{L}_i^M})) = (-1)^{\#S-d+i-1} \operatorname{rank} \widetilde{H}_{\#S-d+i-1}(\Delta([\widehat{0}, S]_{\mathcal{L}_i^M}); \mathbb{Z}).$$

Therefore, for each $1 \leq j \leq m - d + 1 + i$,
$$\begin{aligned}
[p^{d-1-i+j} t^i] \widehat{Y}_M(p,t) &= \omega_j^{i,M} = (-1)^j \sum_{S \in \mathcal{S}_i^M:\ \#S-d+1+i=j} \mu_i^M(S) \\
&= \sum_{S \in \mathcal{S}_i^M:\ \#S-d+1+i=j} \operatorname{rank} \widetilde{H}_{j-2}(\Delta([\widehat{0}, S]_{\mathcal{L}_i^M}); \mathbb{Z}) \\
&= \operatorname{rank} \bigoplus_{S \in \mathcal{S}_i^M} \widetilde{H}_{j-2}(\Delta([\widehat{0}, S]_{\mathcal{L}_i^M}); \mathbb{Z}) \\
&= \operatorname{rank} WH_{j-1}(\mathcal{L}_i^M; \mathbb{Z}),
\end{aligned}$$
completing the proof. □

**Corollary 3.2.** *Let $M = (E, \rho)$ be a matroid of rank $d \geq 1$ with $\#E = m$ and with no loops. For each $0 \leq i \leq d-1$ let $N_i$ be the matroid dual to the $i$-th truncation of $M$. The characteristic polynomial of $M$ is*
$$\chi_M(t) = (-1)^d (1-t) \sum_{i=0}^{d-1} \operatorname{rank} \widetilde{H}_{m-d+i-1}(\mathcal{I}(N_i); \mathbb{Z})(-t)^i.$$

*Proof.* In view of Corollary 2.9 and Theorem 3.1, it suffices to notice that for $0 \leq i \leq d-1$, $WH_{m-d+i}(\mathcal{L}_i^M; \mathbb{Z}) = \widetilde{H}_{m-d+i-1}(\Gamma_i(E); \mathbb{Z})$ (since $E$ is the only set $S \in \mathcal{S}_i^M$ with $\#S - d + 1 + i = m - d + 1 + i$) and that $\Gamma_i(E) = \mathcal{I}(N_i)$. □

Department of Combinatorics and Optimization, University of Waterloo, Waterloo, Ontario, Canada   N2L 3G1

*E-mail address*: dgwagner@math.uwaterloo.ca